\documentclass[11pt]{article}
\usepackage{amssymb,latexsym}
\usepackage{epsfig}
\usepackage{eufrak}
\usepackage{amsmath}
\usepackage{mathrsfs}
\usepackage{color}

\newtheorem{theorem}{Theorem}
\newtheorem{lemma}{Lemma}
\newtheorem{corollary}{Corollary}
\newtheorem{remark}{Remark}

\newcommand{\be}{\begin{equation}}
\newcommand{\ee}{\end{equation}}
\newcommand{\bee}{\begin{eqnarray*}}
\newcommand{\eee}{\end{eqnarray*}}
\newcommand{\bel}{\begin{eqnarray}}
\newcommand{\eel}{\end{eqnarray}}
\newcommand{\bec}{\begin{cases}}
\newcommand{\eec}{\end{cases}}
\newcommand{\bem}{\begin{bmatrix}}
\newcommand{\eem}{\end{bmatrix}}

\newcommand{\bed}{\begin{description}}
\newcommand{\eed}{\end{description}}
\newcommand{\bei}{\begin{itemize}}
\newcommand{\eei}{\end{itemize}}
\newcommand{\ben}{\begin{enumerate}}
\newcommand{\een}{\end{enumerate}}

\newcommand{\beL}{\begin{lemma}}
\newcommand{\eeL}{\end{lemma}}
\newcommand{\beT}{\begin{theorem}}
\newcommand{\eeT}{\end{theorem}}

\newcommand{\bpf}{\begin{pf}}
\newcommand{\epf}{\end{pf}}

\newcommand{\pfbox}{\hfill\mbox{$\Box$}}

\newenvironment{pf}{\paragraph*{Proof{\rm.}}}{\pfbox\bigskip}

\begin{document}

\title{{\bf Order Statistics and Probabilistic Robust
 Control\thanks{This research was supported
in part by grants from AFOSR (F49620-94-1-0415), ARO (DAAH04-96-1-0193), and LE{\cal Q}SF (DOD/LE{\cal Q}SF(1996-99)-04 and LE{\cal Q}SF(1995-98)-RD-A-14).}}}
\author{Xinjia Chen and Kemin Zhou\thanks{All correspondence should be
addressed to this author: Phone: (504) 388-5533, Fax: (504) 388-5200,
email: kemin@ee.lsu.edu}\\
Department of Electrical and Computer Engineering\\
Louisiana State University\\
Baton Rouge, LA 70803\\
chan@ece.lsu.edu \ \ kemin@ee.lsu.edu}

\date{Received in November 12, 1997; Revised in February 25, 1998} \maketitle

\begin{abstract}
  Order statistics theory is applied in this paper
to probabilistic robust control theory to compute the
minimum sample size needed to come up with a reliable estimate
of an uncertain quantity under continuity assumption of
the related probability distribution.  Also, the concept of distribution-free tolerance intervals is applied to estimate the range of an uncertain quantity
and extract the information about its distribution.
To overcome the limitations imposed by the continuity assumption
in the existing  order statistics theory,  we have
derived a cumulative distribution function of the order
statistics without the continuity
assumption and developed an inequality showing that
this distribution has an upper bound which equals to the corresponding distribution  when the continuity assumption is satisfied.
 By applying this inequality, we investigate the minimum computational effort
needed  to come up with an reliable estimate for the upper bound (or lower bound) and the range of a quantity.  We also give
conditions,  which are much weaker than the absolute continuity assumption,
 for the existence of such minimum sample size.
 Furthermore, the issue of making tradeoff
 between performance level and risk is addressed and
 a guideline for making this kind of tradeoff is established. This guideline can be applied in general without continuity assumption.
\end{abstract}


\clearpage

\section{Introduction}
It is now well-known  that many deterministic worst-case
robust analysis and synthesis problems are NP hard, which means that the exact
analysis and synthesis of the corresponding  robust control
problems may be computational demanding \cite{Braatz,Toker,ZDG}.
Moreover,  the deterministic worst-case
robustness measures may be quite conservative due to overbounding of
the system uncertainties.
As pointed out   by
Khargonekar and  Tikku in \cite{KT}, the difficulties of
deterministic worst-case robust control problems
are inherent to the problem formulations and a major change of the
paradigm is necessary. An alternative to the deterministic
approach is the probabilistic approach which has been
studied extensively by Stengel and co-workers, see for example,
 \cite{RS}, \cite{SR}, and references therein.
 Aimed at breaking through the NP-hardness
barrier and reducing the conservativeness of the deterministic
robustness measures, the probabilistic approach  has recently received
a renewed attention in the work by
Barmish and Lagoa \cite{BL},   Barmish, Lagoa, and Tempo \cite{BLT},
Barmish and Polyak \cite{BP}, Khargonekar and  Tikku \cite{KT},
Bai, Tempo, and Fu \cite{bai}, Tempo, Bai,
and Dabbene \cite{TD}, Yoon and Khargonekar \cite{Yoon},
and references therein.

In particular, Tempo, Bai and Dabbene in \cite{TD} and Khargonekar and Tikku in \cite{KT}
have derived bounds for the number of samples
 required to estimate the upper bound of a quantity with
a certain  a priori specified accuracy and confidence.
It is further shown that this probabilistic approach
for robust control analysis and synthesis has low complexity \cite{KT,TD}.
 It should also be pointed out that the uncertain parameters
do  not necessarily have to be random, they can be regarded as
randomized variables as pointed out in \cite{KT}.

One important open question in this area is the minimum computational
effort (i.e., the minimum sample size)
needed to obtain a reliable estimate for
the upper bound (or lower bound, or the range)
of an uncertain quantity.
We shall answer this question in this paper.
This paper is organized as follows. Section 2 presents the
problem formulation and motivations. In Section 3, we derive a
distribution inequality without the continuity assumption.
Design tradeoff is discussed in Section 4 using the distribution
derived in Section 3. Section 5 gives the minimum sample size under various assumptions
and Section 6 considers the tolerance intervals.

 \section{Problem Formulation and Motivations}

Let $q$ be a random vector, bounded in a compact set
${\cal Q}$, with a multivariate probability density function $\varpi(q)$.
Let $u(q)$ be a real scalar measurable function of the random
vector $q$ with cumulative probability distribution $F_{u}(\gamma) := P_{rob}\left\{u(q) \leq \gamma \right\}
$
for a given $\gamma \in {\bf R}$.  Let $q^{1},q^{2},\cdots,q^{N}$ be
the i.i.d. (independent and identically distributed)
samples of $q$ generated according to the same probability
density function $\varpi(q)$ where $N$ is the sample size.
Now define random variable ${\hat{u}}_{i},\; i=1,2,\cdots,N$
as the $i-$th smallest observation of $u(q)$ during $N$ sample experiments.
These random variables are called {\it order statistics} \cite{AB,david}
 because ${\hat{u}}_{1} \leq {\hat{u}}_{2}  \leq {\hat{u}}_{3} \leq \cdots \cdots \leq {\hat{u}}_{N}$.

We are interested in computing the following probabilities:

\begin{itemize}
\item $P_{rob}\left\{P_{rob}\left\{u(q) > {\hat{u}}_{n}\right\} \leq \varepsilon\right\}
$
for any $1 \leq n \leq N$, and $\varepsilon \in(0,1)$;

\item $P_{rob}\left\{P_{rob}\left\{u(q) < {\hat{u}}_{m}\right\} \leq \varepsilon\right\}
$
for any $1 \leq m \leq N$, and $\varepsilon \in(0,1)$;

\item
$P_{rob}\left\{P_{rob}\left\{{\hat{u}}_{m} < u(q)
\leq {\hat{u}}_{n}\right\} \geq 1-\varepsilon \right\}
$
for any $1 \leq m < n \leq N$, and $\varepsilon \in(0,1)$.
\end{itemize}

In the subsequent subsections, we shall give some
motivations for computing such probabilities.

\subsection{Robust Analysis and Optimal Synthesis}

As noted in \cite {KT} and \cite{TD},
to tackle robust analysis or optimal synthesis problem,
it is essential to deal with the following questions:

\begin{itemize}

\item What is $\max_{{\cal Q}}u(q)$
(or $\min_{{\cal Q}}u(q)$)?

\item What is the value of $q$ at which $u(q)$ achieves
 $\min_{\cal Q}u(q)$ (or $\max_{\cal Q}u(q)$)?

\end{itemize}

\begin{figure}[htb] %\label{ff}
\centering{
\setlength{\unitlength}{0.050in}%
\begin{picture}(30,21)
\put(10,10){\framebox(10,10){$M$}}       % M
\put(11,0){\framebox(8,6){$\Delta(q)$}}          % \Delta
\put(5,21){\makebox(0,0){$z$}}           % z
\put(5,8){\makebox(0,0)[r]{$y$}}         % y
\put(25,21){\makebox(0,0){$v$}}          % w
\put(25,8){\makebox(0,0)[l]{$u$}}        % v
\put(10,19){\vector(-1,0){10}}
\put(30,19){\vector(-1,0){10}}
\put(24,13){\vector(-1,0){4}}
\put(6,13){\line(1,0){4}}
\put(19,3){\line(1,0){5}}
\put(6,3){\vector(1,0){5}}
\put(6,3){\line(0,1){10}}
\put(24,3){\line(0,1){10}}
\end{picture}
}
\caption{Uncertain System}
\label{fig_a}
\end{figure}
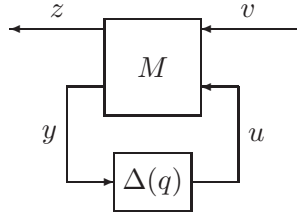

Consider, for example, an uncertain system shown in Figure ~\ref{fig_a}.
Denote the transfer function from $v$ to $z$ by $T_{zv}$ and suppose  that
$T_{zv}$ has the following state space realization
$
T_{zv} = \left[
\begin{array}{c|c}
A (q) & B (q) \\
\hline C(q) & D(q) \end{array} \right].
$
We can now consider following robustness problems:

\begin{itemize}
\item Robust stability: Let
$
u(q) := \max_{i} {\rm Re} \, \lambda_i  (A(q))
$
where $\lambda_i(A)$ denotes the $i$-th eigenvalue of $A$. Then
the system is robustly stable if
$
\max_{q \in {\cal Q}} u(q) < 0.
$
\item Robust performance: Suppose  $A(q)$ is stable for all $q \in {\cal Q}$.
Define
$
u(q):={||T_{zv}||}_{\infty}.
$
Then the robust performance problem is to determine if
$
\max_{q \in {\cal Q}} u(q) \leq \gamma
$
is satisfied for some prespecified $\gamma >0$.

\end{itemize}

\begin{figure}[htb] %\label{fff}
\centering{
\setlength{\unitlength}{0.05in}%
\begin{picture}(30,21)
\put(10,10){\framebox(10,10){$P$}}       % M
\put(11,0){\framebox(8,6){$K(q)$}}          % \Delta
\put(5,21){\makebox(0,0){$z$}}           % z
\put(5,8){\makebox(0,0)[r]{$y$}}         % y
\put(25,21){\makebox(0,0){$d$}}          % w
\put(25,8){\makebox(0,0)[l]{$u$}}        % v
\put(10,19){\vector(-1,0){10}}
\put(30,19){\vector(-1,0){10}}
\put(24,13){\vector(-1,0){4}}
\put(6,13){\line(1,0){4}}
\put(19,3){\line(1,0){5}}
\put(6,3){\vector(1,0){5}}
\put(6,3){\line(0,1){10}}
\put(24,3){\line(0,1){10}}
\end{picture}
}
\caption{Synthesis Framework}
\label{fig_s}
\end{figure}
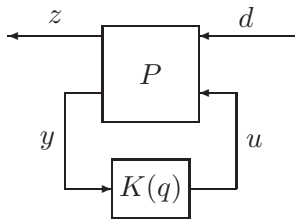

As another example, consider a dynamical system shown in Figure~\ref{fig_s}
and suppose $q$ is a vector of controller parameters to be designed.
Denote the transfer function from $d$ to $z$ by $T_{zd}$
and suppose $T_{zd}$ has the following state space realization
$
T_{zd} = \left[
\begin{array}{c|c}
A_s (q) & B_s (q) \\
\hline C_s(q) & D_s(q) \end{array} \right].
$
  Let
$
u(q):={||T_{zd}||}_{\infty}
$
and
$
{\cal Q}_{s}= \{q \in {\cal Q}:\; A_{s}(q) \; {\rm is \; stable }\}.
$
  Then an $H_\infty$ optimal design problem
is to determine  a vector of parameters $q$ achieving
$
\min_{{\cal Q}_{s}} \; u(q).
$

In general, exactly evaluating
$\min_{\cal Q} \; u(q)$ (or $\max_{\cal Q} \; u(q)$) or
determining $q$ achieving it may be an NP hard problem
and thus is intractable in practice.
Henceforth, we adopt the probabilistic approach
proposed in  \cite{KT} and \cite{TD}.  That is, estimating $\min_{\cal Q} u(q)$
as
\[
{\hat{u}}_{1}=
\min_{ i \in \left\{1,2,\cdots,N \right\}} \;{\hat{u}}_{i}
\]
for sufficiently large $N$ and computing
$P_{rob}\left\{P_{rob}\left\{u(q) < {\hat{u}}_{1}\right\} \leq \varepsilon\right\}
$ for a small $\varepsilon \in (0,1)$ to see how reliable the estimation is.
  Similarly, we estimate $\max_{\cal Q} u(q)$
as
$
{\hat{u}}_{N}=
\max_{ i \in \left\{1,2,\cdots,N \right\} } \;{\hat{u}}_{i}
$
for sufficiently large $N$ and consider $P_{rob}\left\{P_{rob}\left\{u(q) >
 {\hat{u}}_{N}\right\} \leq \varepsilon\right\}
$.

\subsection{Quantity Range}
In many applications, estimating only the upper bound (or lower bound) for a quantity is not sufficient.
It is also important to estimate the range of the quantity with a certain accuracy and confidence level.
For example, in pole placement problem, we need to know the range which the poles fall into.
Suppose that $q$ is the vector of uncertain parameters or design parameters
of a system and $u(q)$ is an uncertain quantity, for example,
$u(q)$ may be the $H_\infty$  norm of a closed-loop transfer function or the
maximum real part of the eigenvalues of the closed-loop system matrix.
Intuitively, the range of quantity $u(q)$ can be approached by $({\hat{u}}_{1}, {\hat{u}}_{N}]$ as sample size $N$ goes to infinity.
  Therefore,
  it is important to know
  $P_{rob}\left\{P_{rob}\left\{{\hat{u}}_{1} < u(q)
  \leq {\hat{u}}_{N}\right\} \geq 1-\varepsilon \right\}$.

So far, we have only concerned the lower bound and (or) upper bound
of uncertain quantity $u(q)$.
Actually, it is  desirable to know its distribution function $F_{u}(.)$.
This is because $F_{u}(.)$
contains all the information of the quantity.  However, the exact computation of
$F_{u}(.)$ is, in general, intractable \cite{XYD}.  An alternative is to
extract as much as possible the information of
the distribution function $F_{u}(.)$
from the observations ${\hat{u}}_{i},\; i=1,2,\cdots,N$. For this purpose,
we are interested in computing the probabilities asked at
the beginning of this section.  In particular, we will see in section $4$ that,
computing $P_{rob}\left\{P_{rob}\left\{u(q) >
{\hat{u}}_{n}\right\} \leq \varepsilon\right\}
$ is of great importance to make the tradeoff between the performance gradation
 and risk when designing a controller.

\section{Distribution Inequality}

Note that $P_{rob}\left\{u(q) \leq {\hat{u}}_{i} \right\} = F_{u}({\hat{u}}_{i}),\;\;i=1,2,\cdots,N
$.  To compute the probabilities asked at the beginning of section $2$, it is important to know
the  associated distribution of any $k$ random variables $F_{u}({\hat{u}}_{{i}_{1}}),\;F_{u}({\hat{u}}_{{i}_{2}}),
\;\cdots,\;F_{u}({\hat{u}}_{{i}_{k}}),
\;\;1 \leq {i}_{1} < {i}_{2} < \cdots < {i}_{k} \leq N,
\;\;1 \leq k \leq N$, i.e.,
\[
F(t_{1},t_{2},\cdots,t_{k})
:= P_{rob}\left\{ F_{u}({\hat{u}}_{i_{1}})
\leq t_{1},\;F_{u}({\hat{u}}_{i_{2}}) \leq t_{2},\;\cdots,\;
F_{u}({\hat{u}}_{i_{k}}) \leq t_{k}\right\}.
\]

 To that end, we have the following theorem which follows essentially by combining {\it Probability Integral Transformation Theorem} in \cite{JG}
and Theorem $2.2.3$ in \cite{david}.
% which is the main result of this paper.
\begin{theorem} \label{main0}
Let $0 = t_{0} \leq t_{1} \leq t_{2} \leq \cdots \leq t_{k} \leq 1$.
Define
\[
f_{ {i}_{1}, {i}_{2}, \cdots , {i}_{k} } (x_{1},x_{2},\cdots,x_{k}):=\prod_{s=0}^{s=k}\;N!\;\frac { {(x_{s+1}-x_{s})}^{i_{s+1}-i_{s}-1} } {(i_{s+1}-i_{s}-1)!}
\]
with $x_{0}:=0$, $x_{k+1}:=1$, $i_{0}:=0$, $i_{k+1}:=N+1$.  Suppose the cumulative distribution function $F_{u}(\gamma) := P_{rob}\left\{u(q) \leq \gamma \right\}
$ is continuous.  Then
\[
F(t_{1},t_{2},\cdots,t_{k}) =
\int_{{\bf D}_{ {t}_{1}, {t}_{2}, \cdots , {t}_{k} }}f_{ {i}_{1}, {i}_{2}, \cdots , {i}_{k} } (x_{1},x_{2},\cdots,x_{k})\;dx_{1} dx_{2} \cdots dx_{k}
\]
where
$
{\bf D}_{ {t}_{1}, {t}_{2}, \cdots , {t}_{k} }:=\left\{(x_{1},x_{2},\cdots,x_{k}):\;\;0 \leq x_{1} \leq x_{2} \leq
\cdots \leq x_{k},\;x_{s} \leq t_{s},\;s=1,2,\cdots,k\right\}.
$
\end{theorem}

\begin{remark}
Theorem~\ref{main0} can play an important role in robust control
as illustrated  in the following sections.
However, its further application is limited by the continuity assumption.
In many robust control problems, it is reasonable to assume that
$u(q)$ is measurable, while the continuity of $F_{u}(\gamma)$ is not
necessarily guaranteed.
For example, $F_{u}(\gamma)$ is not continuous when uncertain quantity $u(q)$ equals to a constant in an open  set of ${\cal Q}$.
We can come up with many uncertain systems with which
the continuity assumption for the distribution of
quantity $u(q)$ is not guaranteed. To tackle these problems
without continuity assumption by probabilistic approach and
investigate the minimum computational effort, we shall
develop a distribution inequality which accommodates
 the case when the continuity is not guaranteed.
\end{remark}

First, we shall established the following lemma.

\begin{lemma} \label{main3} Let U be a random variable with uniform
distribution over $[0,1]$ and ${\hat{U}}_{n},\;\; n=1,2,\cdots,N$ be the order statistics of U,
 i.e., ${\hat{U}}_{1} \leq {\hat{U}}_{2} \leq \cdots \leq {\hat{U}}_{N}$.
Let $0 = t_{0} < t_{1} < t_{2} < \cdots < t_{k} \leq 1$.
Define
$
G_{j_{1},j_{2}, \cdots ,j_{k}}\left(t_{1},t_{2},\cdots,t_{k}\right)
:= (1-t_{k})^{N-\sum_{l=1}^{k}j_{l}}\prod_{s=1}^{k} {N-\sum_{l=1}^{s-1}j_{l} \choose j_{s}}
{(t_{s}-t_{s-1})}^{j_{s}}
$
 and
$
{\bf I}_ { {i}_{1}, {i}_{2}, \cdots , {i}_{k} }:= \left\{(j_{1},j_{2},\cdots,j_{k}):\;\;i_{s} \leq \sum_{l=1}^{s}j_{l} \leq N,\;\;s=1,2,\cdots,k\right\}.
$
Then
\[
P_{rob}\left\{ {\hat{U}}_{i_{1}} \leq t_{1},{\hat{U}}_{i_{2}} \leq t_{2},\cdots,{\hat{U}}_{i_{k}} \leq t_{k}\right\}\\
=\sum_{(j_{1},j_{2},\cdots,j_{k})
\in {\bf I}_ { {i}_{1}, {i}_{2}, \cdots ,
{i}_{k} }}G_{j_{1},j_{2}, \cdots ,j_{k}}\left(t_{1},t_{2},\cdots,t_{k}\right)
\]
\end{lemma}
\begin{pf}
Let $j_{s}$ be the number of samples of U which fall
 into $(t_{s-1},t_{s}],\;\;s=1,2,3,\cdots,k$. Then the number of
 samples of U which fall into $[0,t_{s}]$ is $\sum_{l=1}^{s}j_{l}$.
 It is easy to see that the event
 $\left\{{\hat{U}}_{i_{s}} \leq t_{s}\right\}$
 is equivalent to event $\left\{i_{s} \leq \sum_{l=1}^{s}j_{l} \leq N\right\}$.
   Furthermore,
  the event
  \[\left\{ {\hat{U}}_{i_{1}} \leq t_{1},{\hat{U}}_{i_{2}}
  \leq t_{2},\cdots,{\hat{U}}_{i_{k}} \leq t_{k}\right\}
  \]
  is equivalent to the event
  $
  \left\{i_{s} \leq \sum_{l=1}^{s}j_{l} \leq N,\;\;s=1,2,\cdots,k\right\}.
  $
  Therefore,
\begin{eqnarray*}
&   & P_{rob}\left\{ {\hat{U}}_{i_{1}} \leq t_{1},{\hat{U}}_{i_{2}} \leq t_{2},\cdots,{\hat{U}}_{i_{k}} \leq t_{k}\right\}\\
& = & \sum_{(j_{1},j_{2},\cdots,j_{k}) \in {\bf I}_ { {i}_{1}, {i}_{2}, \cdots , {i}_{k} }}\;
\prod_{s=1}^{k}{N-\sum_{l=1}^{s-1}j_{l} \choose j_{s}}
{(t_{s}-t_{s-1})}^{j_{s}}\;{(1-t_{k})}^{N-\sum_{l=1}^{k}j_{l}}\\
& = & \sum_{(j_{1},j_{2},\cdots,j_{k}) \in {\bf I}_ { {i}_{1}, {i}_{2}, \cdots , {i}_{k} }}\;
G_{j_{1},j_{2}, \cdots ,j_{k}}\left(t_{1},t_{2},\cdots,t_{k}\right).\end{eqnarray*}
\end{pf}

\begin{theorem} \label{main1}
Let $0 = t_{0} \leq t_{1} \leq t_{2} \leq \cdots \leq t_{k} \leq 1$.
Define
\[
\tilde{F}(t_{1},t_{2},\cdots,t_{k}) := P_{rob}\left\{ F_{u}({\hat{u}}_{i_{1}}) < t_{1},\;F_{u}({\hat{u}}_{i_{2}}) < t_{2},\;\cdots,\;F_{u}({\hat{u}}_{i_{k}}) < t_{k}\right\}
\]
and  $\tau_{s}:
=\sup_{ \{x:\;F_{u}(x) <t_{s} \} }F_{u}(x),\;\;s=1,2,\cdots,k$.
Suppose $u(q)$ is a measurable function of $q$. Then
\[
\tilde{F}(t_{1},t_{2},\cdots,t_{k})=
\int_{{\bf D}_{ {\tau}_{1},{\tau}_{2},\cdots,{\tau}_{k} } }
f_{ {i}_{1}, {i}_{2}, \cdots , {i}_{k} } (x_{1},x_{2},\cdots,x_{k})\;dx_{1} dx_{2} \cdots dx_{k}.
\]
Furthermore,
$
\tilde{F}(t_{1},t_{2},\cdots,t_{k}) \leq \int_{{\bf D}_{ {t}_{1}, {t}_{2}, \cdots , {t}_{k} }}f_{ {i}_{1}, {i}_{2}, \cdots , {i}_{k} } (x_{1},x_{2},\cdots,x_{k})\;dx_{1} dx_{2} \cdots dx_{k}
$
and the equality holds if $F_{u}(\gamma)$ is continuous.
\end{theorem}

\begin{pf}
Define $\alpha_{0}:= -\infty$ and
${\alpha}_{s} := \sup\left\{x:\;F_{u}(x)
< t_{s}\right\} ,\;\;{\alpha}_{s}^{-}:= \alpha_{s}-{\epsilon},\;\;s=1,2,\cdots,k$
where $\epsilon > 0$ can be arbitrary small.
Let $\phi_{s}:= F_{u}({\alpha}_{s}^{-}),\;\;s=1,2,\cdots,k$. We can show that
${\phi}_{l} < {\phi}_{s}$ if
${\alpha}_{l} < {\alpha}_{s},\; 1 \leq l < s \leq k$.
  In fact , if this is not true, we have
${\phi}_{l} = {\phi}_{s}$. Because $\epsilon$ can be arbitrary small,
we have $\alpha_{s}^{-} \in (\alpha_{l},\;\alpha_{s})$.  Notice that
${\alpha}_{l} = \min\left\{x:\;F_{u}(x) \geq t_{l} \right\}$
, we have
$t_{l} \leq \phi_{s} = \phi_{l}$.  On the other hand,
by definition we know that
$\alpha_{l}^{-} \in \left\{x:\;F_{u}(x)
< t_{l}\right\}$ and thus $\phi_{l}=F_{u}(\alpha_{l}^{-}) < t_{l}$,
which is a contradiction.  Notice that $F_{u}(\gamma)$ is nondecreasing and right-continuous, we have
${\alpha}_{1} \leq {\alpha}_{2} \leq \cdots \leq {\alpha}_{k}$
and $0 \leq {\phi}_{1} \leq {\phi}_{2} \leq \cdots \leq {\phi}_{k} \leq 1
$
and that event $\left\{F_{u}({\hat{u}}_{i_{s}}) < t_{s}\right\}$ is equivalent to the event $\left\{{\hat{u}}_{i_{s}} < {\alpha}_{s}\right\}$.  Furthermore, event
$
\left\{ F_{u}({\hat{u}}_{i_{1}}) < t_{1},\;F_{u}({\hat{u}}_{i_{2}}) < t_{2},\;\cdots,\;F_{u}({\hat{u}}_{i_{k}}) < t_{k}\right\}
$
is equivalent to event
$\left\{{\hat{u}}_{i_{1}} < {\alpha}_{1},\;{\hat{u}}_{i_{2}} < {\alpha}_{2},\;\cdots,\;{\hat{u}}_{i_{k}} < {\alpha}_{k}\right\}
$
which is defined by $k$ constraints
${\hat{u}}_{i_{s}} < {\alpha}_{s}, \;s=1,2,\cdots,k$.
For every $l <k$, delete constraint ${\hat{u}}_{i_{l}} < {\alpha}_{l}$
if there exists $s >l$ such that ${\alpha}_{s}={\alpha}_{l}$.  Let the remained constraints
be ${\hat{u}}_{i^{'}_{s}} < {\alpha}^{'}_{s}, \;s=1,2,\cdots,k^{'}$ where
${\alpha}^{'}_{1} < {\alpha}^{'}_{2} < \cdots < {\alpha}^{'}_{k^{'}}$.  Since all
constraints deleted are actually redundant, it follows that event $\left\{{\hat{u}}_{i_{1}} < {\alpha}_{1},\;{\hat{u}}_{i_{2}} < {\alpha}_{2},\;\cdots,\;{\hat{u}}_{i_{k}} < {\alpha}_{k}\right\}
$ is equivalent to event
$\left\{{\hat{u}}_{i^{'}_{1}} < {\alpha}^{'}_{1},\;
{\hat{u}}_{i^{'}_{2}} < {\alpha}^{'}_{2},\;\cdots,\;
{\hat{u}}_{i^{'}_{k^{'}}} < {\alpha}^{'}_{k^{'}}\right\}$.  Now let $j_{s}$ be the number of observations of $u(q)$
which fall into $[{\alpha}^{'}_{s-1},{\alpha}^{'}_{s}),\;\;s=1,2,\cdots,k^{'}$.
 Then the number of observations of $u(q)$ which fall into
 $(-\infty,{\alpha}^{'}_{s})$ is $\sum_{l=1}^{s}j_{l}$.
 It is easy to see that the event $\left\{{\hat{u}}_{i^{'}_{s}} <
 {\alpha}^{'}_{s}\right\}$ is equivalent to the event
 $\left\{i^{'}_{s} \leq \sum_{l=1}^{s}j_{l} \leq N\right\}$.
 Furthermore, the event $\left\{{\hat{u}}_{i^{'}_{1}} < {\alpha}^{'}_{1},\;
 {\hat{u}}_{i^{'}_{2}} < {\alpha}^{'}_{2},\;\cdots,\;
 {\hat{u}}_{i^{'}_{k^{'}}} < {\alpha}^{'}_{k^{'}}\right\}$
 is equivalent to event
$
\left\{i^{'}_{s} \leq \sum_{l=1}^{s}j_{l} \leq N,
\;\;s=1,2,\cdots,k^{'}\right\}.
$
Therefore
\begin{eqnarray*}&   & \tilde{F}(t_{1},t_{2},\cdots,t_{k}) =
P_{rob}\left\{ F_{u}({\hat{u}}_{i_{1}}) < t_{1},\;F_{u}({\hat{u}}_{i_{2}}) < t_{2},\;\cdots,\;F_{u}({\hat{u}}_{i_{k}}) < t_{k}\right\}\\
& = & P_{rob}\left\{{\hat{u}}_{i_{1}} < {\alpha}_{1},\;{\hat{u}}_{i_{2}} < {\alpha}_{2},\;\cdots,\;{\hat{u}}_{i_{k}} < {\alpha}_{k}\right\}
= P_{rob}\left\{{\hat{u}}_{i^{'}_{1}} < {\alpha}^{'}_{1},\;
{\hat{u}}_{i^{'}_{2}} < {\alpha}^{'}_{2},\;\cdots,\;
{\hat{u}}_{i^{'}_{k^{'}}} < {\alpha}^{'}_{k^{'}}\right\}\\
& = & \sum_{(j_{1},j_{2},\cdots,j_{k^{'}})
\in {\bf I}_ { {i}^{'}_{1}, {i}^{'}_{2}, \cdots , {i}^{'}_{k^{'}} }}
\prod_{s=1}^{k^{'}}{N-\sum_{l=1}^{s-1}j_{l} \choose j_{s}}
{[F_{u}({{\alpha}^{'}_{s}}^{-})-
F_{u}({{\alpha}^{'}_{s-1}}
^{-})]}^{j_{s}}\;{[1-
F_{u}({{\alpha}^{'}_{k^{'}}}^{-})
]}^{N-\sum_{l=1}^{k^{'}}j_{l}}\\
& = & \sum_{(j_{1},j_{2},\cdots,j_{k^{'}})
\in {\bf I}_ { {i}^{'}_{1}, {i}^{'}_{2}, \cdots , {i}^{'}_{k^{'}} }}
G_{j_{1},j_{2}, \cdots ,j_{k^{'}}}\left(\phi^{'}_{1},
\phi^{'}_{2},\cdots,\phi^{'}_{k^{'}}\right).
\end{eqnarray*}
Now consider event
$\left\{ {\hat{U}}_{i_{1}} \leq \phi_{1},\;{\hat{U}}_{i_{2}}
\leq \phi_{2},\;\cdots,\;{\hat{U}}_{i_{k}} \leq \phi_{k}\right\}$.
For every $l <k$, delete constraint
${\hat{U}}_{i_{l}} \leq {\phi}_{l}$
if there exists $s >l$ such that ${\phi}_{s}={\phi}_{l}$.
Notice that $\phi_{s}=F_{u}({\alpha}_{s}^{-})$ and ${\phi}_{l} < {\phi}_{s}$ if
${\alpha}_{l} < {\alpha}_{s},\; 1 \leq l < s \leq k$, the remained constraints must be
be ${\hat{U}}_{i^{'}_{s}} \leq {\phi}^{'}_{s}, \;s=1,2,\cdots,k^{'}$ where
${\phi}^{'}_{s} = F_{u}({{\alpha}^{'}_{s}}^{-}), \;s=1,2,\cdots,k^{'}
$ and $\phi^{'}_{1} < \phi^{'}_{2} < \cdots < \phi^{'}_{k^{'}}$.  Since all
constraints deleted are actually redundant, it follows that event
$\left\{{\hat{U}}_{i_{1}} \leq {\phi}_{1},\;{\hat{U}}_{i_{2}} \leq
{\phi}_{2},\;\cdots,\;{\hat{U}}_{i_{k}} \leq {\phi}_{k}\right\}
$ is equivalent to event
$\left\{{\hat{U}}_{i^{'}_{1}} \leq {\phi}^{'}_{1},\;
{\hat{U}}_{i^{'}_{2}} \leq {\phi}^{'}_{2},\;\cdots,\;
{\hat{U}}_{i^{'}_{k^{'}}} \leq {\phi}^{'}_{k^{'}}\right\}$.
By Theorem $2.2.3$ in \cite{david} and  Lemma~\ref{main3}
 \begin{eqnarray*}
&   & \int_{{\bf D}_{ {\phi}_{1},{\phi}_{2},\cdots,{\phi}_{k} } }
f_{ {i}_{1}, {i}_{2}, \cdots , {i}_{k} }
(x_{1},x_{2},\cdots,x_{k})\;dx_{1} dx_{2} \cdots dx_{k}\\
& = & P_{rob} \left\{{\hat{U}}_{i_{1}} \leq {\phi}_{1},\;{\hat{U}}_{i_{2}} \leq
{\phi}_{2},\;\cdots,\;{\hat{U}}_{i_{k}} \leq {\phi}_{k}\right\}
= P_{rob}\left\{{\hat{U}}_{i^{'}_{1}} \leq {\phi}^{'}_{1},\;
{\hat{U}}_{i^{'}_{2}} \leq {\phi}^{'}_{2},\;\cdots,\;
{\hat{U}}_{i^{'}_{k^{'}}} \leq {\phi}^{'}_{k^{'}}\right\}\\
& = & \sum_{(j_{1},j_{2},\cdots,j_{k^{'}})
\in {\bf I}_ { {i}^{'}_{1}, {i}^{'}_{2}, \cdots , {i}^{'}_{k^{'}} }}
G_{j_{1},j_{2}, \cdots ,j_{k^{'}}}\left(\phi^{'}_{1},
\phi^{'}_{2},\cdots,\phi^{'}_{k^{'}}\right).
\end{eqnarray*}
Therefore,
$
\tilde{F}(t_{1},t_{2},\cdots,t_{k})=
\int_{{\bf D}_{ {\phi}_{1},{\phi}_{2},\cdots,{\phi}_{k} } }
f_{ {i}_{1}, {i}_{2}, \cdots , {i}_{k} } (x_{1},x_{2},\cdots,x_{k})\;dx_{1} dx_{2} \cdots dx_{k}$.  By the
definitions of $\tau_{s}$ and $\phi_{s}$, we know that
  ${\bf D}_{ {\tau}_{1},{\tau}_{2},\cdots,{\tau}_{k} }$ is the closure of
  ${\bf D}_{ {\phi}_{1}, {\phi}_{2}, \cdots , {\phi}_{k} }$, i.e.,
   $
  {\bf D}_{ {\tau}_{1},{\tau}_{2},\cdots,{\tau}_{k} }
  = {\bar {\bf D}}_{ {\phi}_{1}, {\phi}_{2}, \cdots , {\phi}_{k} }$
  and that their Lebesgue measures
  are equal.  It follows that
 \[
  \tilde{F}(t_{1},t_{2},\cdots,t_{k})
  = \int_{ {\bf D}_{ {\tau}_{1},{\tau}_{2},\cdots,{\tau}_{k} } }
  f_{ {i}_{1}, {i}_{2}, \cdots , {i}_{k} }
  (x_{1},x_{2},\cdots,x_{k})\;dx_{1} dx_{2} \cdots dx_{k}.
  \]
Notice that $\tau_{s} \leq t_{s},\;\;s=1,2,\cdots,k$, we have
  $
  {\bf D}_{ {\tau}_{1},{\tau}_{2},\cdots,{\tau}_{k} }
  \subseteq {\bf D}_{ {t}_{1}, {t}_{2}, \cdots , {t}_{k} }$
  and hence
  \[
  \tilde{F}(t_{1},t_{2},\cdots,t_{k})
  \leq \int_{{\bf D}_{ {t}_{1}, {t}_{2}, \cdots ,
  {t}_{k} }}f_{ {i}_{1}, {i}_{2}, \cdots , {i}_{k} }
  (x_{1},x_{2},\cdots,x_{k})\;dx_{1} dx_{2} \cdots dx_{k}.
  \]
    Furthermore, if $F_{u}(\gamma)$ is continuous,
    then $\tau_{s} = t_{s},\;\;s=1,2,\cdots,k$,
    hence $
  {\bf D}_{ {\tau}_{1},{\tau}_{2},\cdots,{\tau}_{k} }
  = {\bf D}_{ {t}_{1}, {t}_{2}, \cdots , {t}_{k} }$
   and the equality holds.
\end{pf}

\section{Performance and Confidence Tradeoff}

For the synthesis of a controller for an uncertain system, we usually have a conflict between the performance level and robustness.  The following theorem helps to make the tradeoff.
\begin{theorem} \label{ccce}
Let $1 \leq n \leq N,\;1 \leq m \;\leq N$ and $\varepsilon \in(0,1)$. Suppose $u(q)$ is measurable.
 Then
\be \label{ine1}
P_{rob}\left\{P_{rob}\left\{u(q) > {\hat{u}}_{n}\right\} \leq \varepsilon\right\} \geq 1- \int_{0}^{1-\varepsilon} \frac{N!}{(n-1)!(N-n)!}x^{n-1}(1-x)^{N-n}dx
\ee
and  the equality holds if and only if
$\sup_{ \{x:\;F_{u}(x) < 1- \varepsilon \} } F_{u}(x) = 1-\varepsilon$;
Moreover,
\be  \label{ine2}
P_{rob}\left\{P_{rob}\left\{u(q) < {\hat{u}}_{m}\right\}
\leq \varepsilon\right\} \geq 1-
\int_{0}^{1-\varepsilon} \frac{N!}{(m-1)!(N-m)!}x^{N-m}(1-x)^{m-1}dx
\ee
and the equality holds if and only if
$\inf_{\{ x:\;F_{u}(x) > \varepsilon \}} F_{u}(x) = \varepsilon$.
\end{theorem}
\begin{pf}
Apply Theorem ~\ref{main1} to the case of $k=1,\;i_{1}=n$, we have
\begin{eqnarray*}
&   & P_{rob}\left\{F_{u}({\hat{u}}_{n}) < 1-\varepsilon\right\}\\
& = & \int_{{\bf D}_{\tau}} \frac{N!}{(n-1)!(N-n)!}x^{n-1}(1-x)^{N-n}dx \leq \int_{0}^{1-\varepsilon} \frac{N!}{(n-1)!(N-n)!}x^{n-1}(1-x)^{N-n}dx
\end{eqnarray*}
where ${\bf D_{\tau}}=(0, \tau]$
with $\tau=\sup_{ \{ x:\;F_{u}(x) < 1-\varepsilon \} } F_{u}(x)$.
Therefore,
\begin{eqnarray*}
&   & P_{rob}\left\{P_{rob}\left\{u(q) > {\hat{u}}_{n}\right\} \leq \varepsilon\right\}=P_{rob}\left\{F_{u}({\hat{u}}_{n}) \geq 1-\varepsilon\right\}\\
& = & 1-P_{rob}\left\{F_{u}({\hat{u}}_{n}) < 1-\varepsilon\right\} \geq 1-\int_{0}^{1-\varepsilon} \frac{N!}{(n-1)!(N-n)!}x^{n-1}(1-x)^{N-n}dx.
\end{eqnarray*}
The equality holds if and only if
$\sup_{ \{x:\;F_{u}(x) < 1- \varepsilon \} } F_{u}(x) = 1-\varepsilon$ because ${\bf D_{\tau}}=(0, 1-\varepsilon]$ if and only if $\tau=\sup_{ \{ x:\;F_{u}(x) < 1-\varepsilon \} } F_{u}(x) = 1-\varepsilon$.

Now let $v(q)=-u(q)$.  Let the cumulative distribution function of $v(q)$ be
$F_{v}(.)$ and define order statistics ${\hat{v}}_{i},\;\;i=1,2\cdots,N$ as the
$i$-th smallest observation of $v(q)$ during $N$ i.i.d. sample experiments,
i.e.,
${\hat{v}}_{1} \leq {\hat{v}}_{2}  \leq {\hat{v}}_{3} \leq \cdots \cdots
\leq {\hat{v}}_{N}$.  Obviously,
 ${\hat{u}}_{m}=-{\hat{v}}_{N+1-m}$ for any $1 \leq m \leq N$.
It is also clear that $F_{v}(-x)=1- F_{u}(x^{-})$,
which leads to the equivalence of
$
\sup_{ \{x:\;F_{v}(x) < 1-\varepsilon\} }F_{v}(x)=1-\varepsilon
$
and
$
\inf_{\{ x:\;F_{u}(x) > \varepsilon \}}F_{u}(x)=\varepsilon$.  Therefore,
apply ~(\ref{ine1}) to the situation of $v(q)$, we have
\begin{eqnarray*}
&   & P_{rob}\left\{P_{rob}\left\{u(q) < {\hat{u}}_{m}\right\}
\leq \varepsilon\right\}\\
& = & P_{rob}\left\{P_{rob}\left\{v(q) > {\hat{v}}_{N+1-m}\right\}
\leq \varepsilon\right\} \geq 1-
\int_{0}^{1-\varepsilon} \frac{N!}{(m-1)!(N-m)!}x^{N-m}(1-x)^{m-1}dx
\end{eqnarray*}
and the equality holds if and only if
$\inf_{\{ x:\;F_{u}(x) > \varepsilon \}} F_{u}(x) = \varepsilon$.
\end{pf}

\begin{figure}[htb]
\centerline{\psfig{figure=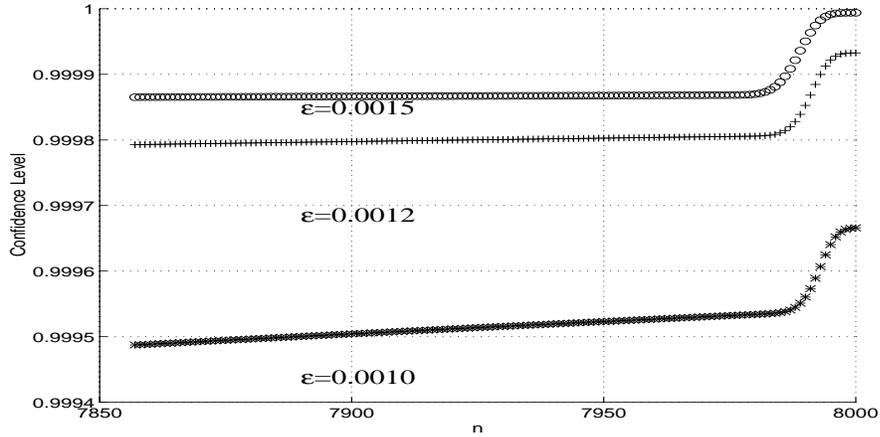,height=2.3in,width=4.6in}}
\caption{Confidence Level $:=P_{rob}\left\{P_{rob}\left\{u(q) >
{\hat{u}}_{n} \right\} \leq \varepsilon\right\}$.  Sample size
$N=8000$.  Performance level increases as $n$ decreases, while
confidence level decreases.} \label{fig_bound22}
\end{figure}

In Figure \ref{fig_bound22} we computed the lower bound for $P_{rob}\left\{P_{rob}\left\{u(q) > {\hat{u}}_{n} \right\} \leq \varepsilon\right\}
$ for sample size $N=8000$ with $\varepsilon = 0.0010$, $\varepsilon = 0.0012$, and $\varepsilon=0.0015$
respectively.

Theorem~\ref{ccce} can be  used as a guideline to robust control analysis and synthesis. For example, when dealing with
 robust stability problem, we need to compute the maximum of the
 real part of the closed-loop poles, denoted by $u(q)$, which is a function of uncertain parameters $q$. When we estimate the upper bound of $u(q)$ by sampling, it is possible that most of the samples concentrated in an interval and very few fall far beyond that interval. If we take ${\hat{u}}_{N}$ as the maximum in the design of a controller, it may be conservative. However, if we choose $n$ such that ${\hat{u}}_{n}$ is much smaller than  ${\hat{u}}_{N}$, while
$P_{rob}\left\{P_{rob}\left\{u(q) > {\hat{u}}_{n} \right\} \leq \varepsilon\right\}
$
is close to $1$,
then the controller based on  ${\hat{u}}_{n}$ may have
much better performance but with only a little bit more increase of risk.
For example, let's say, sample size $N=8000$ and $\varepsilon=0.0010$, the distribution of sample is like this, ${\hat{u}}_{1},\;{\hat{u}}_{2},\;\cdots,\;{\hat{u}}_{7900}$ concentrated in an interval, and ${\hat{u}}_{7900}$ is much smaller than ${\hat{u}}_{7901},\;{\hat{u}}_{7902},\;\cdots,{\hat{u}}_{8000}$. It is sure that the controller designed by taking ${\hat{u}}_{7900}$ as the upper bound will have much higher performance level than the controller designed by taking ${\hat{u}}_{8000}$ as the upper bound. To compare the risks for these two cases, we have
$P_{rob}\left\{P_{rob}\left\{u(q) > {\hat{u}}_{8000} \right\} \leq 0.0010\right\} \geq 0.99966
$
and $P_{rob}\left\{P_{rob}\left\{u(q) > {\hat{u}}_{7900} \right\} \leq 0.0010\right\} \geq 0.99951.
$
These data indicate that there is only a little bit more increase of risk
by taking ${\hat{u}}_{7900}$ instead of ${\hat{u}}_{8000}$ as
 the upper bound in designing a controller.

\section{Minimum Sample Size}
In addition to
the situation of making the tradeoff between performance degradation
 and risk, Theorem~\ref{ccce} can also play an important role in the issue of
computational effort required to come up with an estimate of
the upper bound (or lower bound) of a quantity with a certain accuracy and confidence.
This issue was first addressed independently  by Tempo, Bai and Dabbene in \cite{TD} and Khargonekar
 and Tikku in \cite{KT} and their results are summarized in
the following theorem.
\begin{theorem} \label{cccc}
For any $\varepsilon,\delta\in(0,1)$,
$P_{rob}\left\{P_{rob}\left\{u(q) >
{\hat{u}}_{N}\right\} \leq \varepsilon\right\} \geq  1-\delta$ if
$N \geq \frac{\ln{\frac{1}{\delta}}}{\ln{\frac{1}{1-\varepsilon}}}$.
\end{theorem}
\begin{remark}
This theorem only answers the question that how much computational effort is sufficient.
We shall also concern about what is the minimum computational effort.
By applying Theorem ~\ref{ccce} to the case of $n=N$,
we can also obtain Theorem ~\ref{cccc}.  Moreover, we can see that,
for a certain accuracy (i.e., a fixed value of $\varepsilon$), the bound becomes minimum if and only if
\be
\sup_{ \{ x:\;F_{u}(x) < 1-\varepsilon \} }F_{u}(x)=1-\varepsilon.
\label{mar}
\ee
If $F(u)$ is continuous (i.e., ~(\ref{mar}) is guaranteed for
any $\varepsilon \in (0,1)$), the bound is of course tight.

Similarly, by Theorem ~\ref{ccce} to the case of $m=1$ we have that
for $\varepsilon,\delta\in(0,1)$,
\[
P_{rob}\left\{P_{rob}\left\{u(q) \geq {\hat{u}}_{1}\right\}
\geq 1-\varepsilon\right\} \geq  1-\delta
\]
if $N \geq \frac{\ln{\frac{1}{\delta}}}{\ln{\frac{1}{1-\varepsilon}}}$.
For a fixed $\varepsilon \in(0,1)$, this bound is tight if and only if
 $
  \inf_{\{ x:\;F_{u}(x) > \varepsilon \}} F_{u}(x)=\varepsilon.
  $
\end{remark}

\section{Quantity Range and Distribution-Free Tolerance Intervals}
 To estimate the range of an uncertain quantity with a certain accuracy and confidence level apriori specified, we have the following corollary.

\begin{corollary} \label{ccc}
Suppose $F_{u}(\gamma)$ is continuous.  For any $\varepsilon,\delta\in(0,1)$,
\[
P_{rob}\left\{P_{rob}\left\{{\hat{u}}_{1} < u(q) \leq {\hat{u}}_{N}\right\} \geq 1-\varepsilon\right\} \geq  1-\delta
\]
if and only if $\mu(N) \leq  \delta$ where $\mu(N) := {(1-\varepsilon)}^{N-1}\left[1+(N-1)\varepsilon\right]$ is a monotonically decreasing function of $N$.
\end{corollary}

\begin{remark} The minimum $N$ guaranteeing this condition can be found
 by a simple {\it bisection} search. This bound of sample size is minimum because our computation of probability is exact.
This bound is also practically small, for example,
$N \geq 1,483$ if $\varepsilon = \delta = 0.005$,
and  $N \geq 9,230$ if $\varepsilon = \delta = 0.001$.
Therefore, to obtain a reliable estimate of the range of an uncertain quantity, computational complexity is not an issue.
\end{remark}

In general, it is important to know the probability of a quantity falling between two arbitrary samples. To that end, we have

\begin{corollary} \label{general}
Suppose $F_{u}(\gamma)$ is continuous.  Then for $\varepsilon \in(0,1)$ and $1 \leq m < n \leq N$,
\[
P_{rob}\left\{P_{rob}\left\{{\hat{u}}_{m} < u(q)
\leq {\hat{u}}_{n}\right\} \geq 1-\varepsilon \right\}=
\int_{1-\varepsilon}^{1} N {N-1 \choose n-m-1}
 x^{n-m-1}{(1-x)}^{N-n+m} dx.
\]
\end{corollary}
Here $({\hat{u}}_{m},\;{\hat{u}}_{n}]$ is referred as distribution-free tolerance interval in order statistics theory (see \cite{david}).

\end{document}